\documentclass[11pt]{article}
\usepackage{amsmath}
\usepackage{amssymb}
\usepackage{latexsym}

\usepackage{amsthm, amstext}
\usepackage{array, amsfonts, mathrsfs}
\usepackage{graphicx}

\usepackage{enumitem}

\usepackage[numbers,sort&compress]{natbib}\bibpunct[, ]{[}{]}{,}{n}{,}{,}
\usepackage[english]{babel}  
\usepackage[utf8]{inputenc}

\theoremstyle{plain}
\newtheorem{theorem}{Theorem}

\begin{document}

\title{On the Cauchy problem for the wave equation in a two-dimensional domain with data on the boundary}

\author{M.~N.~Demchenko\footnote{St.~Petersburg Department of
V.\,A.~Steklov Institute of Mathematics of
the Russian Academy of Sciences, 
27 Fontanka, St.~Petersburg, Russia. E-mail: demchenko@pdmi.ras.ru.\newline
\indent The research was supported by the RFBR grant 20-01-00627-a.}}

\date{}

\maketitle
\begin{flushright}
{\em\large
\begin{tabular}{r}
Dedicated to Vasiliy Mikhaylovich Babich\\ 
on the occasion of his jubilee
\end{tabular}
}
\end{flushright}

\begin{abstract}
The subject of the paper is the Cauchy problem for 
the wave equation in a space-time cylinder $\Omega\times{\mathbb R}$,
$\Omega\subset{\mathbb R}^2$, with the data on the surface
$\partial\Omega\times I$, where $I$ is a finite time interval.
The algorithm for solving the Cauchy problem
with the data 
on 
$S\times I$, $S\subset\partial\Omega$, was obtained previously.
Here we adapt
this algorithm to the special case
$S=\partial\Omega$ and show that 
in this situation, the solution is determined
with higher stability in comparison with the case 
$S\subsetneqq\partial\Omega$.

\smallskip

\end{abstract}

\section{Introduction}
We deal with the problem of determination of a function $u(x,y,t)$
satisfying the wave equation
\begin{equation}
  \partial_t^2 u - \Delta u = 0 \label{wave}
\end{equation}
in the cylinder $(x,y,t)\in\Omega\times{\mathbb R}$
($\Omega$ is a bounded domain in ${\mathbb R}^2$)
from the Cauchy data on a surface that is 
a subset of
$\partial\Omega\times{\mathbb R}$.
In contrast to the classical Cauchy problem with data given
at $t=t_0$, the problem in consideration is ill-posed.
However, 
the solution $u$
is determined uniquely, at least in a subset of the space-time cylinder~\cite{LRSh, Tat, is}.
Our goal is to construct an algorithm for solving the Cauchy problem in consideration.
One of such algorithms, which applies in the case when $\Omega$
is the half-space of arbitrary dimension, 
was first proposed by R.~Courant (see~\cite{K}).
In~\cite{Blag}, 
the inversion formula for the problem in the 
three-dimensional half-space was obtained, which
allows determining a solution to the wave equation
from the Cauchy data given on a certain
unbounded subset of the space-time boundary.
The case of an arbitrary Euclidean domain was considered 
in~\cite{DD18, mn99}.
The algorithm obtained there allows determining the solution 
$u$ from the Cauchy data given on the set
$S\times I$, where $S$ is a part of the boundary $\partial\Omega$
($S\subsetneqq\partial\Omega$), $I$ is a finite interval.
The solution $u$ is determined in some subset of the cylinder
$\Omega\times{\mathbb R}$ dependent on $S$ and $I$.
Note also that there are 
results on determination
of singularities of the solution from the Cauchy data 
on the boundary, see~\cite{Quintoetal} and the literature
cited therein.
Typical tools for this are the ray method (see~\cite{BB, BK}),
Maslov's canonical operator, Fourier integral operators.
Problems of determination of solutions to hyperbolic equations from boundary data arise in
geophysics~\cite{Kab} 
and coefficient inverse problems~\cite{bel}.

To construct an algorithm of determining the solution $u$,
we will use the results of papers~\cite{DD18, mn99} mentioned above,
where the Cauchy problem with the data on a part of the boundary of the domain was considered.
We will 
adapt the algorithm obtained there to the specific case
when the data are given on the whole boundary
(in the notation introduced above, this means that $S=\partial\Omega$).
We will show that in this situation, the algorithm is more stable
than it is in the case $S\subsetneqq\partial\Omega$.
Note that in the three-dimensional case,
the Cauchy data on $\partial\Omega\times I$
are sufficient for the stable reconstruction of the solution
at a given time by Kirchhoff's formula (see, e.g., \cite{FPR}).
The latter is derived by applying Green's formula to
the solution $u$ and the fundamental solution of the wave equation.
However, in the two-dimensional problem, 
which is considered in the present paper, 
this approach requires the Cauchy data on an unbounded
time interval.
The reason for this is the fact that 
the corresponding fundamental solution
(unlike that in the three-dimensional case)
is non-zero
for arbitrarily large times for any fixed $x$, $y$.

We will consider the case when the solution $u$
satisfies the Dirichlet boundary condition
\begin{equation}
  u|_{\partial\Omega\times{\mathbb R}} = 0.
  \label{Dir}
\end{equation}
Thus the Cauchy data reduce to the values of the derivative $\partial_\nu u$
of the solution along the outward unit normal to the boundary
of the domain.

\section{On the Cauchy problem with the data on a part 
of the boun\-dary}
Here we outline the approach proposed in~\cite{DD18, mn99}
and give some formulas from these papers,
which will be used further.
In this section, the Dirichlet boundary condition~(\ref{Dir})
is not required.

The algorithm obtained in~\cite{DD18, mn99} is based
on a special choice of a function $V$ in 
Green's formula 
\begin{equation}
  \int_{\Omega\times{\mathbb R}} [u\, (\partial^2_t - \Delta) V - 
    V\, (\partial^2_t - \Delta) u]\, dx dy dt =
  \int_{\partial\Omega\times{\mathbb R}} \left[V\, \partial_\nu u 
    - u\, \partial_\nu V \right] d\sigma dt.
  \label{Green}
\end{equation}
Introduce the function, which depends on the parameter $h>0$,
\begin{equation}
  v_h(x,y,t) = -\theta(y-|t|) w_h(x,y,t).
  \label{tew}
\end{equation}
Here $\theta$ is the Heaviside function,
and $w_h(x,y,t)$ is defined for $y\geqslant|t|$ as
\begin{align}
  w_h(x,y,t) = \frac{1}{4\pi^{3/2}\sqrt{h}} \int_{-\pi}^{\pi} 
  \exp\left(-\frac{1}{h}\left(x-i\sqrt{y^2-t^2}\,
  \cdot {\sin}s\right)^2\right)\, ds \label{vrep}
\end{align}
(here we give the formula only for the two-dimensional case).
The following relation holds true~\cite{DD18, mn99}
\begin{align}
  &\partial^2_t v_h - \Delta v_h =
  \frac{e^{-x^2/h}}{\sqrt{\pi h}} \delta(y)\delta(t).
  \label{wavedelta}
\end{align}
Note that the first factor on the right hand side
is the Gaussian distribution on the line,
which converges to $\delta(x)$ as $h\to 0$.
If we formally substitute $v_h$ for $V$ in Green's formula~(\ref{Green}),
we obtain
\begin{align*}
      &\int_{\Omega\times{\mathbb R}} u\, (\partial^2_t - \Delta) v_h\, dx dy dt
  =\int_{\partial\Omega\times{\mathbb R}} \left[v_h \partial_\nu u 
    - u\, \partial_\nu v_h\right] d\sigma dt.
\end{align*}
The derivatives of the function $v_h$ are distributions.
Here we do not specify the meaning of the integrals
of these distributions in the given equality;
this issue is discussed in the cited papers.
In view of~(\ref{wavedelta}), sending $h\to 0$ yields
\begin{align*}
  &u(0,0,0) = 
  \lim_{h\to 0} 
  \int_{\partial\Omega\times{\mathbb R}} \left[v_h \partial_\nu u 
    - u\, \partial_\nu v_h\right] d\sigma dt.
\end{align*}
Next we make use of the fact that the function
$v_h(x,y,t)$ and its derivatives converge to zero
as $h\to 0$, if the point $(x,y)$ lies outside the cone
$\{ y \geqslant |x| \}$.
This can be derived by analyzing the exponent on the right hand side
of~(\ref{vrep}).
Taking into account also that for fixed
$x$, $y$, 
the function $v_h$ vanishes for sufficiently large $t$,
we conclude that on the right hand side of the relation
obtained above,
the set of integration $\partial\Omega\times{\mathbb R}$ can be replaced by
$S\times I$, where 
$I\subset{\mathbb R}$ is a sufficiently large time interval,
and $S\subset\partial\Omega$ is an arbitrary relatively open set
that contains the intersection of the specified cone with $\partial\Omega$.
Thus we obtain the following equality
\begin{equation}
  u(0,0,0) = 
  \lim_{h\to 0} 
  \int_{S\times I} \left[v_h \partial_\nu u 
    - u\, \partial_\nu v_h\right] d\sigma dt,
  \label{local}
\end{equation}
which allows solving the Cauchy problem with the data on a part of the boundary.

The Cauchy problem can be solved by other choices of
$V$ in Green's formula~(\ref{Green}).
The simplest modification of the algorithm described above
can be obtained by passing from $v_h$ to the function
\begin{equation}
  v_h^\alpha(x,y,t) = v_h(x',y',t),
  \label{tewa}
\end{equation}
where the transformation
\begin{equation}
  (x,y) \mapsto (x',y') = 
  (x\cos\alpha + y\sin\alpha, -x\sin\alpha + y\cos\alpha)
  \label{transcoord}
\end{equation}
is the rotation through the angle $\alpha$ around the origin.
This yields the analogue of formula~(\ref{local}),
in which the set $S$ should contain all of the points
$(x,y)\in\partial\Omega$ such that $y' \geqslant |x'|$.
In particular, this means that the solution $u$ at the same point
can be determined from the Cauchy data given on various sets.
This is a typical feature of ill-posed Cauchy problems
(including problems for elliptic equations).

Now we make some remarks on the stability of 
the algorithm~(\ref{local}).
On the set of integration $S\times I$, the function $v_h$
grows exponentially as $h\to 0$, which can be deduced from 
formulas~(\ref{tew}), (\ref{vrep}).
Therefore, if the Cauchy data $u$, $\partial_\nu u$ are given
with some error, the limit on the right hand side of~(\ref{local})
can be considerably different from the value of the solution or may not exist at all.
Thus the right hand side should be approximated by the value
of the integral for some positive $h$.
This integral is close to the ``convolution'' of the solution
with the right hand side of~(\ref{wavedelta}).
The latter has a positive dispersion (in $x$),
which tends to zero as $h\to 0$ and determines
the approximation accuracy. 
Thus one can increase the approximation accuracy by decreasing $h$,
which, on the other hand,
decreases the stability of the algorithm 
(in case when the Cauchy data are given with some error)
because of the growth of the function $v_h$.

To solve the Cauchy problem with the data on the whole boundary
$\partial\Omega$, rather than on a part of it,
we will substitute a certain function $V_h$ 
to Green's formula,
which grows 
not as rapidly as $v_h$ does when $h$ goes to zero.
This provides a more stable reconstruction of the solution.
The corresponding dispersion (see formulas~(\ref{rhode}))
is close to that of the right hand side of~(\ref{wavedelta}) 
for $v_h$, which means that the accuracy
is principally the same as in formula~(\ref{local}).
Our construction of the function $V_h$ relies
on the functions $v_h$, $w_h$. 

\section{The function $V_h$}\label{choiceV}
Define the function $V_h$ as 
the average of the function $v_h^\alpha$, the latter being defined 
by relation~(\ref{tewa}), over $-\pi\leqslant\alpha\leqslant\pi$
\begin{equation}
  V_h(x,y,t) = \frac{1}{2\pi}\int_{-\pi}^\pi v_h^\alpha(x,y,t)\, d\alpha.
  \label{Vdef}
\end{equation}
Clearly the function $V_h$ is radially symmetric with respect to $x$, $y$.
The equality
\begin{equation}
  V_h(x,y,t) = 0, \quad x^2+y^2 < t^2
  \label{Vsupp}
\end{equation}
follows from the analogous equality for $v_h^\alpha$,
which is valid for any $\alpha$. 
Due to~(\ref{wavedelta}), we have
\[
  (\partial^2_t v_h^\alpha - \Delta v_h^\alpha)(x,y,t) = 
  \frac{e^{-x'^2/h}}{\sqrt{\pi h}} \delta(y')\delta(t).
\]
Integrating this equality over $-\pi\leqslant\alpha\leqslant\pi$, 
we find
\begin{equation}
  \partial_t^2 V_h - \Delta V_h = \rho_h(x,y) \delta(t),
  \quad
  \rho_h(x,y) = \frac{1}{\pi \sqrt{\pi h}}\,
  \frac{e^{-(x^2+y^2)/h}}{\sqrt{x^2+y^2}}.
  \label{rhode}
\end{equation}
The function $\rho_h$ has a singularity at the origin.
However, this singularity is integrable
and we have
\begin{equation}
  \rho_h(x,y)\to \delta(x)\delta(y), \quad h\to 0
  \label{derho}
\end{equation}
in the sense of distributions, which can be verified 
by a direct calculation.
Therefore the right hand side of the first relation in~(\ref{rhode}) converges to
$\delta(x)\delta(y)\delta(t)$ as $h\to 0$.
Recall that due to~(\ref{wavedelta}), the function $v_h$
also enjoys this property, which was used 
in the derivation of~(\ref{local}).

\begin{theorem}\label{t}
  Suppose that $\Omega$ is a bounded domain in ${\mathbb R}^2$, $\partial\Omega\in C^\infty$, and the function $u\in C^2(\overline\Omega\times{\mathbb R})$ satisfies
  the wave equation~(\ref{wave}) and the boundary condition~(\ref{Dir}). 
  Then for any $(x^*,y^*)\in\Omega$, $t^*\in{\mathbb R}$, we have
  \begin{align}
    u(x^*,y^*,t^*) = \lim_{h\to 0} 
    \int_{\partial\Omega} d\sigma_{x,y} \int_{t^*-\tau(x,y)}^{t^*+\tau(x,y)} 
    V_h^*\, \partial_\nu u\, dt,
    \label{result}
  \end{align}
  where $V_h^*(x,y,t) = V_h(x-x^*,y-y^*,t-t^*)$,
  $\tau(x,y)$ is the distance between
  $(x,y)$ and $(x^*,y^*)$ in the Euclidean metric.
\end{theorem}

In addition to giving a proof of Theorem~\ref{t},
we will also investigate the behavior of the function $V_h$
as $h\to 0$.
We will show that 
$V_h$ 
satisfies the following estimate
\begin{equation}
  V_h(x,y,t) = O(h^{-1/2}), \quad h\to 0
  \label{Vest}
\end{equation}
locally uniformly in $x$, $y$, $t$,
which is in contrast to the exponential growth of $v_h$.
In the course of the derivation of this estimate,
we will obtain a representation for $V_h$, which,
unlike~(\ref{Vdef}),
does not involve expressions that grow exponentially fast
as $h\to 0$.

\section{A modified representation for $w_h$}
Our definition of $V_h$ involves the functions $v_h$, $w_h$.
The latter is defined for $y\geqslant |t|$ by formula~(\ref{vrep}).
Clearly this formula makes sense also when $|y|\geqslant |t|$
and, moreover, when $|y| < |t|$.
It suffices to observe that in the latter case,
the choice of the branch 
of the square root in the exponent 
in~(\ref{vrep}) is immaterial since the factor $\sin s$
changes sign on the set of integration
$-\pi\leqslant s\leqslant\pi$.
Hence for $|y|<|t|$, we may write
\begin{align}
  w_h(x,y,t) = \frac{1}{4\pi^{3/2}\sqrt{h}} \int_{-\pi}^{\pi} 
  \exp\left(-\frac{1}{h}\left(x+\sqrt{t^2-y^2}\,
  \cdot {\sin}s\right)^2\right)\, ds, \label{vrepout}
\end{align}
assuming the square root in the exponent to be positive.

We will need a certain modification of relations~(\ref{vrep}), (\ref{vrepout}). 
Suppose that $y,\, t\geqslant 0$.
In the case $y > t$, the substitution 
$\xi = e^{i s} \sqrt{(y+t)/|y-t|}$ in the integral on the right hand side of~(\ref{vrep}) yields the following integral
\begin{align}
  \int 
  \exp\left(-\frac{1}{h}\left(x - 
  \frac{(y-t)\xi}{2} + \frac{y+t}{2\xi}\right)^2\right)
  \frac{d\xi}{i\xi}.
  \label{wxi}
\end{align}
For $y<t$, the same substitution in the integral in~(\ref{vrepout})
yields
\begin{align*}
  \int 
  \exp\left(-\frac{1}{h}\left(x + 
  \frac{i (y-t)\xi}{2} + \frac{i (y+t)}{2\xi}\right)^2\right) 
  \frac{d\xi}{i\xi}.
\end{align*}
Both integrals are taken over the counterclockwise oriented circle
\[
  |\xi|=\sqrt{(y+t)/|y-t|}.
\]
Note that the last integral takes the form~(\ref{wxi}) after
the substitution $\xi' = -i \xi$,
and so in both cases $y>t$ and $y<t$,
it is sufficient to deal with the expression~(\ref{wxi}).
After the deformation of the path of integration to the circle
$|\xi| = 1$, we obtain the following relation
\[
  w_h(x,y,t) = \frac{1}{4\pi^{3/2}\sqrt{h}}
  \int_{|\xi|=1}
  \exp\left(-\frac{1}{h}\left(x - 
  \frac{(y-t)\xi}{2} + \frac{y+t}{2\xi}\right)^2\right)
  \frac{d\xi}{i\xi}
\]
being valid for all $y,\,t\geqslant 0$,
including the case $y=t$, which is treated by the continuity
with respect to $y$, $t$.
Now make the substitution $\xi=e^{i s}$, $-\pi\leqslant s \leqslant \pi$, 
in the integral on the right hand side. We have
\[
  x - \frac{(y-t)\xi}{2} + \frac{y+t}{2\xi} = 
  x - \frac{y}{2}\left(\xi-\frac{1}{\xi}\right) +
  \frac{t}{2}\left(\xi+\frac{1}{\xi}\right) = 
  x - i y \sin s + t \cos s,
\]
whence
\begin{align}
  w_h(x,y,t) = \frac{1}{4\pi^{3/2}\sqrt{h}} \int_{-\pi}^\pi 
  \exp\left(-\frac{1}{h}\left(x - i y \sin s + t \cos s\right)^2\right) 
  ds. \label{intds}
\end{align}
In the resulting representation for the function $w_h$,
the restriction $y,\,t\geqslant 0$ can be dropped since both
sides 
are even with respect to each of the 
variables $x$, $y$, $t$.

\section{Estimate of the function $V_h$}
In this section, we investigate the behavior of 
the function $V_h$ as $h\to 0$.
Our goal is to obtain estimate~(\ref{Vest}).
Since the function $V_h$ is radially symmetric with respect to 
$x$, $y$,
we may assume $x=0$, $y>0$.
Moreover, in view of~(\ref{Vsupp}), we may also assume $y>|t|$.
According to~(\ref{transcoord}) we have
\[
  x' = y\sin\alpha, \quad y' = y\cos\alpha.
\]
Now~(\ref{tew}) implies that $v_h^\alpha(0,y,t) = 0$ for $y'<|t|$,
that is, for $y\cos\alpha < |t|$. 
Hence, setting $\alpha_0 = \arccos(|t|/y)$, we obtain
\begin{align}
  V_h(0,y,t) = \frac{1}{2\pi} 
  \int_{-\alpha_0}^{\alpha_0} v_h^\alpha(0,y,t)\, d\alpha =
  -\frac{1}{2\pi} 
  \int_{-\alpha_0}^{\alpha_0} w_h(y\sin\alpha, y\cos\alpha, t)\, d\alpha.
  \label{Vw}
\end{align}
It is convenient to extend the set of integration 
to the interval $-\pi/2 \leqslant \alpha \leqslant \pi/2$,
so we set
\begin{align*}
  \widetilde V_h(y,t) = -\frac{1}{2\pi} 
  \int_{-\pi/2}^{\pi/2} w_h(y\sin\alpha, y\cos\alpha, t)\, d\alpha.
\end{align*}
Applying representation~(\ref{intds}) for $w_h$, we obtain
\begin{align}
  \widetilde V_h(y,t) = 
-\frac{c}{\sqrt{h}}
  \int_{-\pi}^\pi ds
  \int_{-\pi/2}^{\pi/2}
  \exp\left(-\frac{1}{h}\left(y\,(\sin\alpha - i \sin s \cos\alpha) 
  + t \cos s\right)^2\right) d\alpha,
  \label{V0}
\end{align}
where $c=1/(8\pi^{5/2})$.

\begin{figure}[b!]\centering
  \includegraphics[width=.5\textwidth]{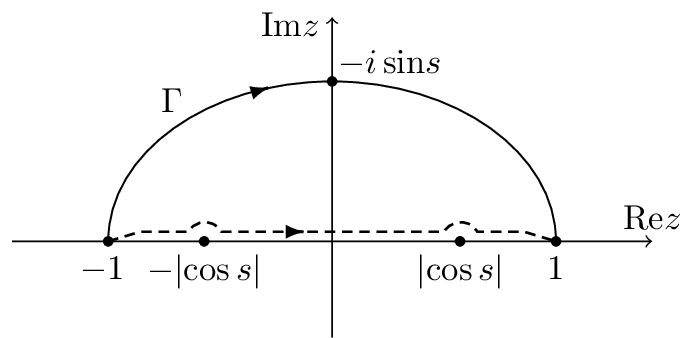}
\caption{The path $\Gamma$ in the case
  $\sin s < 0$ (the case $\sin s > 0$ requires an obvious modification). The dashed line shows the deformed path that bypasses the branch points $\pm\cos s$.}
  \label{figcont}
\end{figure}

Consider the inner integral in~(\ref{V0}) in the case $s \ne \pm\pi/2$.
Make the substitution $z = \sin\alpha - i \sin s \cos\alpha$.
For $-\pi/2 \leqslant \alpha \leqslant \pi/2$, the variable $z$
ranges over the path $\Gamma$ shown on the Fig.~\ref{figcont}.
We have
\[
  \sin\alpha = \frac{z + i \sin s \,
    \sqrt{(\cos s)^2-z^2}}{(\cos s)^2},
  \quad
  \cos\alpha = \frac{-i z \sin s + 
    \sqrt{(\cos s)^2-z^2}}{(\cos s)^2}.
\]
Note that in case if the branch points $z = \pm \cos s$ 
belong to the path $\Gamma$ (which occurs only if $|\cos s| = 1$),
they coincide with its ends.
Now substitution $\alpha=0$ to the last formulas
ensures
we should choose the principal branch of the square root
(for arbitrary $s$).
We have
\[
  dz = (\cos\alpha + i \sin s\, \sin\alpha)\, d\alpha = 
  \sqrt{(\cos s)^2-z^2}\, d\alpha,
\]
Therefore the integral with respect to $\alpha$ in~(\ref{V0}) 
equals 
\begin{align*}
  \int_\Gamma
  \frac{\exp\left(-(y z + t \cos s)^2/h\right) dz}{\sqrt{(\cos s)^2-z^2}}.
\end{align*}
Next we make a deformation of the path $\Gamma$
to $[-1\pm i 0, 1\pm i 0]$, where the sign is chosen to be 
opposite to that of $\sin s$ (in the case $\sin s=0$,
the choice of the sign is immaterial since the branch points
$\pm\cos s = \pm 1$ coincide with the ends of the path).
Denote temporarily $E = \exp\left(-(y z + t \cos s)^2/h\right)$.
For $\sin s < 0$, we have
\begin{align*}
  \int_{-1+i 0}^{1+i 0}
  \frac{E dz}{\sqrt{(\cos s)^2-z^2}} =& 
  \int_{-|\cos s|}^{|\cos s|}
  \frac{E dz}{\sqrt{(\cos s)^2-z^2}}
  + \int_{|\cos s|}^1
  \frac{E dz}
       {-i \sqrt{z^2-(\cos s)^2}}\notag\\
  &+ \int_{-1}^{-|\cos s|}
  \frac{E dz}
       {i \sqrt{z^2-(\cos s)^2}}.
\end{align*}
In the case $\sin s > 0$, the corresponding integral over
$z\in[-1- i\, 0, 1- i\, 0]$ 
is given by a similar expression,
which differs from the previous one in the signs
before the last two terms.
Thus in the integral over $-\pi \leqslant s \leqslant \pi$ in~(\ref{V0}),
these two terms for any $s$ are canceled out by the corresponding
terms for $-s$.
Now we may write equality~(\ref{V0}) in the following form
\begin{align}
  \widetilde V_h(y,t) =& 
  -\frac{c}{\sqrt{h}}
  \int_{-\pi}^\pi ds
  \int_{-|\cos s|}^{|\cos s|}
  \frac{E dz}{\sqrt{(\cos s)^2-z^2}} \notag\\
  &=
  -\frac{4 c}{\sqrt{h}}
  \int_0^{\pi/2} ds
  \int_{-\cos s}^{\cos s}
  \frac{\exp\left(-(y z + t \cos s)^2/h\right) dz}{\sqrt{(\cos s)^2-z^2}}.
  \label{V0fin}
\end{align}

Representation~(\ref{V0fin}) shows that the function $\widetilde V_h$
equals $O(h^{-1/2})$ as $h\to 0$.
According to~(\ref{Vw}) and the definition of $\widetilde V_h$,
for the difference between $V_h(0,y,t)$ and $\widetilde V_h(y,t)$
we have
\[
  \widetilde V_h(y,t) - V_h(0,y,t) = 
  -\frac{1}{\pi} \int_{\alpha_0}^{\pi/2} 
  w_h(y\sin\alpha, y\cos\alpha, t)\, d\alpha
\]
(we used the fact that the function $w_h$ is even in each of the variables).
Since $y\cos\alpha < |t|$ on the set of integration,
representation~(\ref{vrepout}) for $w_h$ is valid,
in which, as in~(\ref{V0fin}), 
the exponent is nonpositive. 
Thus the last relation together with formula~(\ref{V0fin}) for
$\widetilde V_h$ yield a representation for the function $V_h$,
which does not contain expressions that grow exponentially fast
as $h\to 0$
(recall that in the original representation~(\ref{Vdef}),
the function $v_h^\alpha$ grows exponentially fast).
This implies, in particular, the validity of estimate~(\ref{Vest}).

\section{The function $V_h^\varepsilon$}\label{hom}
From representation~(\ref{Vw}) one can deduce that
the function $V_h$ is locally bounded in ${\mathbb R}^2\times{\mathbb R}$,
continuous at $(x,y,t)\ne 0$, and smooth
outside the surface $\{x^2+y^2=t^2\}$.
However, the first derivatives of $V_h$ are unbounded in
the vicinity of the specified surface,
and $V_h$ is discontinuous
at $(x,y,t) = 0$.
Since in the proof of Theorem~\ref{t},
Green's formula will be used, we will pass to smoother
convolutions 
\[
  V_h^\varepsilon(x,y,t) = (\psi^\varepsilon * V_h)(x,y,t) =
  \int_{\mathbb R} \psi^\varepsilon(t-t_1) V_h(x,y,t_1) dt_1
\]
(here and further $*$ denotes the convolution with respect to $t$)
with a mollifier
\[
  \psi^\varepsilon(t) = \frac{e^{-t^2/\varepsilon}}{\sqrt{\pi\varepsilon}}, \quad \varepsilon>0.
\]
In view of~(\ref{Vsupp}), this convolution is well defined, and we have
\[
  V_h^\varepsilon(x,y,t) = 
  \int_{-\tau(x,y)}^{\tau(x,y)} \psi^\varepsilon(t-t_1) V_h(x,y,t_1) dt_1,
\]
where $\tau(x,y) = \sqrt{x^2+y^2}$ coincides with the function
$\tau$ defined in the formulation of Theorem~\ref{t}
providing that $x^*=y^*=0$.
It follows from this representation that the function
$V_h^\varepsilon(x,y,t)$ is continuous in ${\mathbb R}^2\times{\mathbb R}$.
It is also smooth with respect to $t$ for any fixed $(x,y)$.
The derivatives
$\partial_t^m V_h^\varepsilon = (\partial_t^m\psi^\varepsilon) * V_h$, $m\geqslant 0$,
which can be represented similarly to $V_h^\varepsilon$, 
are continuous in ${\mathbb R}^2\times{\mathbb R}$ as well.
Moreover, for any compact set $K\subset{\mathbb R}^2$, we have
\begin{equation}
  |\partial_t^m V_h^\varepsilon(x,y,t)| \leqslant C(K,m,\varepsilon,h), \quad (x,y,t)\in K\times{\mathbb R}.
  \label{bound}
\end{equation}

Next we establish that
$V_h^\varepsilon(\cdot, \cdot, t_0) \in C^\infty({\mathbb R}^2\setminus\{0\})$
for any $t_0\in{\mathbb R}$.
For a function $\eta(x,y)$ from $C_0^\infty({\mathbb R}^2)$
define $\tilde\eta(x,y,t)$ as follows
\[
  \tilde\eta(x,y,t) = \psi^\varepsilon(t_0-t) \eta(x,y).
\]
We have
\begin{align*}
  &\langle{}V_h, (\partial_t^2 - \Delta) \tilde\eta\rangle = 
  \int_{{\mathbb R}^2} dx dy\, \eta(x,y) 
  \int_{\mathbb R} (\partial_t^2\psi^\varepsilon)(t_0-t)\, V_h(x,y,t) \, dt \\
  & - \int_{{\mathbb R}^2} dx dy\, \Delta\eta(x,y) 
  \int_{\mathbb R} \psi^\varepsilon(t_0-t)\, V_h(x,y,t) \, dt \\
  & = \langle(\partial_t^2 V_h^\varepsilon)(\cdot,\cdot,t_0), \eta\rangle -
  \langle{}V_h^\varepsilon(\cdot, \cdot, t_0), \Delta\eta\rangle.
\end{align*}
Here and further angle brackets mean 
pairing of a distribution with a test function.
Note that $\tilde\eta$, unlike $\eta$,
is not a compactly supported function 
due to the factor $\psi^\varepsilon(t_0-t)$, which is nonzero
for all $t$.
However, the pairing of $V_h$ with $\tilde\eta$ and its derivatives
makes sense as the intersection of the supports
${\rm supp} V_h$ and ${\rm supp}\, \tilde\eta$ is compact.
The latter follows from~(\ref{Vsupp}) and the fact that 
the support of $\eta$ is compact.
Equality~(\ref{rhode}) yields
\[
  \langle{}V_h, (\partial_t^2 - \Delta) \tilde\eta\rangle =
  \langle\rho_h, \tilde\eta(\cdot, \cdot, 0)\rangle =
  \psi^\varepsilon(t_0) \langle\rho_h, \eta\rangle.
\]
Combining this with the previous equality, we find
\[
  \langle{}V_h^\varepsilon(\cdot, \cdot, t_0), \Delta\eta\rangle =
  \langle(\partial_t^2 V_h^\varepsilon)(\cdot,\cdot,t_0), \eta\rangle -
  \psi^\varepsilon(t_0) \langle\rho_h, \eta\rangle.
\]
Thus we establish
the following relation, understood in the generalized sense
\begin{equation}
  \Delta V_h^\varepsilon(x, y, t_0) = 
  (\partial_t^2 V_h^\varepsilon)(x,y,t_0) - 
  \psi^\varepsilon(t_0) \rho_h(x,y).
  \label{DeltaV}
\end{equation}
The first term on the right hand side can be estimated 
using~(\ref{bound}).
The second term contains the function $\rho_h$,
which, according to~(\ref{rhode}), has a singularity at the origin.
However, this function is locally bounded in ${\mathbb R}^2\setminus\{0\}$.
Hence~(\ref{DeltaV}) means that for any compact set
$K\subset{\mathbb R}^2\setminus\{0\}$, we have
\begin{equation}
  |\Delta V_h^\varepsilon(x,y,t)| \leqslant C(K,\varepsilon,h), \quad (x,y,t)\in K\times{\mathbb R}.
  \label{boundDelta}
\end{equation}

Taking into account the equality $\psi^\varepsilon = \psi^{\varepsilon/2} * \psi^{\varepsilon/2}$,
we have
\[
  \partial_{t}^2 V_h^\varepsilon = 
  (\partial_t^2 \psi^\varepsilon) * V_h = 
  (\partial_t^2 \psi^{\varepsilon/2}) * V_h^{\varepsilon/2}.
\]
Thus we may write the generalized Laplace operator
of the right hand side of~(\ref{DeltaV}) as follows
\[
  \Delta^2 V_h^\varepsilon(x, y, t_0) = 
  ((\partial_t^2 \psi^{\varepsilon/2}) * \Delta V_h^{\varepsilon/2})(x,y,t_0) - 
  \psi^\varepsilon(t_0) \Delta \rho_h(x,y).
\]
Again the right hand side can be estimated 
at $(x,y)\in K\subset{\mathbb R}^2\setminus\{0\}$.
For this, we apply estimate~(\ref{boundDelta}),
substituting $\varepsilon$ with $\varepsilon/2$.
This results in the analogous estimate for $\Delta^2 V_h^\varepsilon$.
Iterating this argument, we obtain this estimate for
any order of the Laplace operator
applied to $V_h^\varepsilon$.
According to the theory of elliptic equations~\cite{Lad, Evans},
this implies that 
$V_h^\varepsilon(\cdot,\cdot,t_0) \in C^\infty({\mathbb R}^2\setminus\{0\})$.

In addition to the smoothness of $V_h^\varepsilon(\cdot,\cdot,t_0)$ outside the origin, we have also established that
the function $\Delta V_h^\varepsilon(\cdot,\cdot,t_0)$
is locally integrable in ${\mathbb R}^2$, which is a consequence
of equality~(\ref{DeltaV}), estimate~(\ref{bound}), 
and formula~(\ref{rhode}) for $\rho_h$.
This allows us to justify Green's formula applied to $V_h^\varepsilon$.
Let a bounded domain $\Omega$ with a smooth boundary
contain the origin, 
and let $\eta$ be an arbitrary function from $C^2(\overline\Omega)$.
The following equality holds true
\begin{equation}
  \int_\Omega [V_h^\varepsilon\, \Delta \eta - (\Delta V_h^\varepsilon)\, \eta]\,
  dx dy =
  \int_{\partial\Omega} [V_h^\varepsilon \partial_\nu \eta - \eta \partial_\nu V_h^\varepsilon]\,d\sigma.
  \label{GreenDelta}
\end{equation}
(here and further in this argument,
$V_h^\varepsilon$ means $V_h^\varepsilon(\cdot,\cdot,t_0)$).
On the right hand side, we have an integral of a smooth function.
Due to what was said above,
the integral on the left hand side is absolutely convergent.
To establish~(\ref{GreenDelta}), we suppose first that
$\eta\in C^\infty(\overline\Omega)$.
Let $\chi_0$ be a $C^\infty$-smooth function 
with a compact support in $\Omega$,
such that $\chi_0=1$ in a neighborhood of the origin.
Put $\chi_1 = 1 - \chi_0$, and write
$\eta = \chi_0\eta +\chi_1\eta$.
Since $\chi_0\eta$ is a smooth compactly supported function
in ${\mathbb R}^2$,
by definition of the generalized Laplace operator we have
\[
  \int_\Omega V_h^\varepsilon\, \Delta (\chi_0 \eta)\, dx dy =
  \int_{{\mathbb R}^2} V_h^\varepsilon\, \Delta (\chi_0 \eta)\, dx dy =
  \int_{{\mathbb R}^2} (\Delta V_h^\varepsilon)\, \chi_0 \eta\, dx dy.
\]
The function $V_h^\varepsilon$ is smooth on ${\rm supp}\chi_1$,
so we may integrate by parts, which yields
\[
  \int_\Omega V_h^\varepsilon\, \Delta (\chi_1 \eta)\, dx dy =
  \int_{\partial\Omega} [V_h^\varepsilon \partial_\nu \eta - \eta \partial_\nu V_h^\varepsilon]\,d\sigma
  + \int_\Omega (\Delta V_h^\varepsilon)\, \chi_1 \eta\, dx dy
\]
(in the integral over $\partial\Omega$, we took into account that
$\chi_1|_{\partial\Omega} = 1$).
Summing these two equalities gives us relation~(\ref{GreenDelta}),
which is generalized to $C^2$-smooth functions $\eta$
by approximation.

\section{Proof of Theorem~\ref{t}}\label{proof}
With no loss of generality we may assume that the domain $\Omega$
contains the origin, and $x^*=y^*=t^*=0$.

Applying Green's formula~(\ref{GreenDelta}) to the function
$\eta = u(\cdot,\cdot,t)$
and integrating the resulting equality over $-L<t<L$
($L$ is an arbitrary positive number for the time being)
yields
\[
  \int_{-L}^L dt \int_\Omega [V_h^\varepsilon\, \Delta u - 
    u\, \Delta V_h^\varepsilon]\, dx dy =
  \int_{-L}^L dt \int_{\partial\Omega} V_h^\varepsilon\, \partial_\nu u\, d\sigma.
\]
We took into account the boundary condition~(\ref{Dir}).
Since $V_h^\varepsilon$ is smooth with respect to $t$ for any 
fixed $x$, $y$, we have
\begin{align*}
  \int_{-L}^L dt \int_\Omega [u\, \partial^2_t V_h^\varepsilon - 
    V_h^\varepsilon\, \partial^2_t u]\,dx dy
  &= \sum_\pm \int_\Omega (\pm 1) 
  [u \partial_t V_h^\varepsilon - V_h^\varepsilon\, \partial_t u]\big|_{t=\pm L}\, dx dy.
\end{align*}
This equality relies on the fact that it is possible to
change the order of integration with respect to
$t$ and $x$, $y$.
The latter is 
due to the fact that $V_h^\varepsilon$
and all its derivatives with respect to $t$ 
are continuous in ${\mathbb R}^2\times{\mathbb R}$,
as indicated in sec.~\ref{hom}.
The last two equalities imply
\begin{align*}
  &\int_{-L}^L dt \int_\Omega [u\, (\partial^2_t - \Delta) V_h^\varepsilon - 
    V_h^\varepsilon\, (\partial^2_t - \Delta) u]\, dx dy \notag\\
  &=
  \int_{-L}^L dt \int_{\partial\Omega} V_h^\varepsilon\, \partial_\nu u\, d\sigma +
  \sum_\pm \int_\Omega (\pm 1) [u\, \partial_t V_h^\varepsilon - V_h^\varepsilon \partial_t u]\big|_{t=\pm L}\, dx dy.
\end{align*}
Further we denote the second term on the right hand side by $R$.
Taking into account~(\ref{DeltaV})
and the wave equation~(\ref{wave}), 
the last relation takes the following form
\begin{equation}
  \int_{-L}^L dt\, \psi^\varepsilon(t) \int_\Omega u(x,y,t)\, \rho_h(x,y)\, dx dy 
  = 
  \int_{-L}^L dt \int_{\partial\Omega} V_h^\varepsilon\, \partial_\nu u \, d\sigma + R.
  \label{Green2}
\end{equation}
Now we send $\varepsilon$ to zero. 
It follows from~(\ref{Vsupp}) that
$V_h^\varepsilon(x,y,t) \to 0$ as $\varepsilon\to 0$, 
if $x^2+y^2 < t^2$. 
This is also true for $\partial_t V_h^\varepsilon(x,y,t)$ since
$\partial_t V_h^\varepsilon(x,y,t) = \partial_t \psi^\varepsilon * V_h$.
Pick a constant $L$ such that $L>{\rm diam}\,\Omega$.
Then $R\to 0$ as $\varepsilon\to 0$.

Consider the first term on the right hand side of~(\ref{Green2}).
As indicated in the beginning of sec.~\ref{hom},
the function $V_h$ is continuous at $(x,y,t)\ne 0$.
Hence on the set of integration,
the convolutions $V_h^\varepsilon$ are uniformly bounded 
with respect to $\varepsilon$ and converge to $V_h$ pointwise as $\varepsilon\to 0$.
This implies the same assertion for the integrand
$V_h^\varepsilon\, \partial_\nu u$, so the integral tends to
\[
  \int_{-L}^L dt \int_{\partial\Omega} V_h\, \partial_\nu u \, d\sigma =
  \int_{\partial\Omega} d\sigma_{x,y} \int_{-\tau(x,y)}^{\tau(x,y)} 
    V_h\, \partial_\nu u\, dt
\]
as $\varepsilon\to 0$.
In the last equality we used~(\ref{Vsupp}).
The function $\tau$ defined in the formulation of Theorem~\ref{t}
equals $\sqrt{x^2+y^2}$, which is due to our assumption $x^*=y^*=0$. 
Similar argument applies to the left hand side of~(\ref{Green2}).
Therefore passing to the limit $\varepsilon\to 0$ in equality~(\ref{Green2}) yields
\[
  \int_{\Omega} u(x,y,0)\, \rho_h(x,y)\, dx dy = 
  \int_{\partial\Omega} d\sigma_{x,y} \int_{-\tau(x,y)}^{\tau(x,y)} 
    V_h\, \partial_\nu u\, dt.
\]
In view of~(\ref{derho}),
passing to the limit $h\to 0$ in the resulting relation
leads to identity~(\ref{result}).

\end{document}